\documentclass[11pt]{amsart}
\usepackage{amssymb, url}
\usepackage[all]{xy}
\newtheorem{thm}{Theorem} 
\newtheorem{prop}[thm]{Proposition}
\newtheorem{lem}[thm]{Lemma}

\theoremstyle{definition}

\newtheorem{remark}[thm]{Remark}

\newtheorem{example}[thm]{Example}

\newcommand{\spec}{\operatorname{Spec}}
\newcommand{\Char}{\operatorname{char}}

\newcommand{\cO}{\mathcal{O}_x(X)}
\newcommand{\ccO}{\widehat{\mathcal{O}}_x(X)}

\newcommand{\GL}{\operatorname{GL}}
\newcommand{\SL}{\operatorname{SL}}

\newcommand{\Aut}{\operatorname{Aut}}

\newcommand{\id}{\operatorname{id}}
\newcommand{\Sp}{\operatorname{Sp}}
\newcommand{\Gal}{\operatorname{Gal}}

\newcommand{\trdeg}{\operatorname{trdeg}}

\newcommand{\diag}{\operatorname{diag}}
\newcommand{\lra}{\longrightarrow}

\newcommand{\Mat}{{\operatorname{M}}}

\newcommand{\PGLn}{{\operatorname{PGL}}_n}
\newcommand{\PGL}{{\operatorname{PGL}}}
\newcommand{\Sym}{\operatorname{S}}

\newcommand{\cF}{\mathcal{F}}

\newcommand{\bbZ}{\mathbb{Z}}

\newcommand{\bbP}{\mathbb{P}}
\newcommand{\bbA}{\mathbb{A}}
\newcommand\cE{\mathcal{E}}
\newcommand{\SO}{\operatorname{SO}}
\newcommand{\Orth}{\operatorname{O}}
\newcommand{\Br}{\operatorname{Br}}
\newcommand{\inn}{\operatorname{inn}}

\begin{document}

\title[Birational isomorphisms between twisted group actions]{Birational 
isomorphisms between twisted group actions}

\author{Zinovy Reichstein}
\address{Department of Mathematics, University of British Columbia,
Vancouver, BC, Canada, V6T 1Z2}
\email{reichst@math.ubc.ca}
\urladdr{www.math.unbc.ca/~reichst}
\thanks{Z. Reichstein was supported in part by an NSERC research grant}
\thanks{A. Vistoli was partially supported by the University of Bologna, funds for selected research topics}

\author{Angelo Vistoli}
\address{Department of Mathematics, Dipartimento di Matematica,
Universita` di Bologna,
Piazza di Porta S.Donato 5, 40137 Bologna, Italy}
\email{vistoli@dm.unibo.it}
\urladdr{www.dm.unibo.it/~vistoli}



\begin{abstract}
Let $X$ be an algebraic variety with a generically free action
of a connected algebraic group $G$. Given an automorphism 
$\phi \colon G \lra G$, we will denote by $X^{\phi}$ 
the same variety $X$ with the $G$-action given by
$g \colon  x \lra \phi(g) \cdot x$.

V. L. Popov asked if $X$ and $X^{\phi}$ are always $G$-equivariantly
birationally isomorphic.
We construct examples to show that this is not the case in general.
The problem of whether or not such examples can exist in the case 
where $X$ is a vector space with a generically free linear 
action, remains open.
On the other hand, we prove that $X$ and $X^{\phi}$ are always stably 
birationally isomorphic, i.e., $X \times \bbA^m$ and 
$X^{\phi} \times \bbA^m$ are $G$-equivariantly birationally 
isomorphic for a suitable $m \ge 0$.
\end{abstract}

\maketitle
\tableofcontents

\section{Introduction}

Throughout this note 
all algebraic varieties, algebraic groups, group actions and maps 
between them will be defined over a fixed base field $k$.
By a $G$-variety $X$ we shall mean an algebraic variety with a (regular)
action of a linear algebraic group $G$. A morphism (respectively,
rational map, birational isomorphism, etc.) of $G$-varieties is
a $G$-equivariant morphism (respectively, rational map, birational 
isomorphism).  

Given an automorphism $\phi$ of $G$, we can ``twist" a
group action $\alpha \colon G \times X \lra X$ by $\phi$ to
obtain a new $G$-action $\alpha^{\phi}$ on $X$ as follows: 
\[ \alpha^{\phi} \colon G \times X 
\stackrel{(\phi, id)}{\lra} G \times X \stackrel{\alpha}{\lra} X \, . \]
Note that the new action has the same orbits as the old one.
If $X$ is a $G$-variety (via $\alpha$) then we will denote the
``twisted" $G$-variety (i.e., $X$ with the action given by $\alpha^{\phi}$)
by $X^{\phi}$.

Now suppose that the action $\alpha$ of $G$ on $X$ is generically 
free, i.e., that there exists a $G$-invariant open dense 
subset $U$ of $X$ such that the stabilizer of every geometric 
point of $U$ is trivial.  V. L. Popov asked if $X$ and $X^{\phi}$ 
are always birationally isomorphic as $G$-varieties.  
Of particular interest to him was the case where $X$ is a linear 
representation of $G$.  In this case Popov's question is closely 
related to a conjecture of Katsylo~\cite{katsylo}, which says 
that generically free linear $G$-representations $V$ and $W$ 
are ($G$-equivariantly) birationally 
isomorphic if and only if $\dim(V) = \dim(W)$. Katsylo's 
conjecture is known to be false for some finite groups $G$; 
in particular, there are counterexamples, where 
$W = V^{\phi}$ for an automorphism $\phi$ of $G$; 
see~\cite[Section 7]{ry}. On the other hand, 
Katsylo's conjecture remains open for many finite 
groups (e.g., for the symmetric groups $G = \Sym_n$, 
$n \ge 5$) and for all connected semisimple groups.

Two simple observations are now in order.
First note that only the class of $\phi$ in 
the group of outer automorphisms of $G$ matters here.
Indeed, suppose $\phi' = \phi \circ \inn_h$, 
where $\inn_h \colon G \lra G$
is conjugation by $h \in G$, i.e., 
$\phi(g) = h g h^{-1}$. Then 
$X^{\phi}$ and $X^{\phi'}$ are isomorphic via
$x \mapsto h \cdot x$. 
In particular, $X$ and $X^{\phi}$ are
always isomorphic if $G$ has no outer automorphisms, e.g.,
if $G$ is the full symmetric group $\Sym_n$ ($n \ne 6$) 
or if $G$ is semisimple algebraic group whose 
Dynkin diagram has no non-trivial automorphisms; 
cf.~\cite[Theorem 25.16]{boi}. The latter class of groups 
includes every (almost) simple algebraic group, other than those 
of type $A_n$, $D_n$ and $E_6$; cf.~\cite[pp. 354 -- 355]{boi}.

Secondly, the $G$-invariant rational functions for 
$X$ and $X^{\phi}$ are exactly the same, i.e., 
\begin{equation} \label{e.invariants}
k(X)^G = k(X^{\phi})^G \subset k(X) \, . 
\end{equation}
Recall that the inclusion $k(X) \subset k(X)^G$ induces a dominant
rational map $\pi \colon X \dasharrow X/G$, which is called
{\em the rational quotient map}. The $k$-variety $X/G$ is defined
(up to birational isomorphism) by the $k(X/G) = k(X)^G$. 
Thus~\eqref{e.invariants} can be rephrased by saying that
a rational quotient map $\pi$ for $X$ is also a rational
quotient map for $X^{\phi}$.  Note that by a theorem of 
Rosenlicht~\cite{rosenlicht1},~\cite{rosenlicht2}, 
$\pi^{-1}(x)$ is a single $G$-orbit for $x \in X/G$ 
in general position; cf. also~\cite[Section 2.4]{pv}.

The main results of this note are Theorems~\ref{thm1} and~\ref{thm2}
below.

\begin{thm} \label{thm1} Let $X$ be a generically free $G$-variety 
and $\phi \colon G \lra G$ be an automorphism of $G$. Then 
the $G$-varieties $X$ and $X^{\phi}$ are stably birationally 
isomorphic. More precisely there exists an integer $m \ge 0$
and a birational isomorphism of 
\[ f \colon X \times \bbA^m \dasharrow X^{\phi} \times \bbA^m \]
such that the diagram 
\[ \xymatrix{
  X \times \bbA^m \ar@{-->}[r]^{f} \ar@{-->}[d] &
                X^{\phi} \times \bbA^m \ar@{-->}[d]  \cr
  X/G \ar@{=}[r] & X^{\phi}/G,}
\]
commutes. 
\end{thm}

Here $\bbA^m$ is the $m$-dimensional affine space with trivial 
$G$-action, and the vertical map on the left is the composition
of the projection $X \times \bbA^m \lra X$ with the rational 
quotient map $X \dasharrow X/G$ (similarly for the vertical 
map on the right).

\smallskip
\begin{thm} \label{thm2} Let $n \ge 3$ and let
$\phi$ be the (outer) automorphism of $G = \PGLn$ given by 
$g \lra (g^{-1})^{\rm transpose}$. Assume the base field $k$ 
is infinite and contains a primitive $n$th root of unity.  Then 
there exists a generically free $\PGLn$-variety $X$ such that
$X$ and $X^{\phi}$ are not birationally isomorphic over $k$.
\end{thm}

The problem of whether or not such examples can exist 
in the case where $X$ is a vector space with a generically
free linear action of a connected linear algebraic group $G$,
remains open.

\section{The no-name lemma}
\label{sect.no-name}

Recall that $G$-bundle $\pi \colon E \lra X$ is an algebraic vector bundle
with a $G$-action on $E$ and $X$ such that $\pi$ is $G$-equivariant
and the action of every $g \in G$ restricts to a linear
map $\pi^{-1}(x) \lra \pi^{-1}(gx)$ for every $x \in X$. 

Our proof of Theorem~\ref{thm1} in the next section will
heavily rely the following result.

  \begin{lem}[No-name Lemma]
  \label{no-name}
  Let $\pi \colon E \lra X$ be a $G$-bundle of rank
  $r$. Assume that the $G$-action on $X$ is generically free.
  Then there exists a birational isomorphism
  $\pi \colon E \stackrel{\simeq}{\dasharrow}
  X \times \bbA^r$ of $G$-varieties such that the following diagram commutes
  \begin{equation} \label{e.no-name1}
  \xymatrix{ E \ar@{-->}[r]^{\phi \quad }
  \ar@{->}[d]^{\pi} & X \times \bbA^r
  \ar@{->}[dl]^{\operatorname{pr}_1}  \cr X & }
  \end{equation}
  Here $G$ is assumed to act trivially on $\bbA^r$, and
  $\operatorname{pr}_1$ denotes the projection to the first factor.
  \end{lem}

The term ``no-name lemma" is due 
to Dolgachev~\cite{dolgachev}.
In the case where $G$ is a finite group, a proof 
can be found, e.g., in~\cite[Proposition 1.1]{em}, 
\cite[Proposition 1.3]{lenstra} or \cite[Appendix 3]{shafarevich}. 
For a proof in the case where 
the base field $k$ is algebraically closed, $\Char(k) = 0$, and $G$ is 
an arbitrary linear algebraic group, 
see~\cite{bk}, \cite{katsylo},~\cite[Section 4]{cgr}.

In the sequel we would like to use Lemma~\ref{no-name}
in the case where $k$ is not necessarily algebraically closed.
With this in mind, we will prove a more general variant 
of this result (Proposition~\ref{no-name1} below).
For the rest of this section we will work over an arbitrary 
base field $k$. 

\begin{remark}\label{restriction} 
Suppose $G$ is a group scheme of finite type over $k$, 
$X$ is an arbitrary quasi-separated scheme (or algebraic space) 
over $k$ on which $G$ acts (quasi-separated means that the diagonal 
embedding $X \hookrightarrow X \times_{\spec k}X$ is quasi-compact;\
this is automatically satisfied when $X$ is of finite type over $k$). 
We say that the action is \emph{free} when the stabilizers 
of all geometric point of $X$ 
are trivial (as group schemes). This is equivalent to saying that 
the morphism $G \times_{\spec k} X \lra X \times_{\spec k} X$ defined 
in functorial terms by $(g,x) \mapsto (gx, x)$ is categorically 
injective (equivalently, it is injective on geometric points 
and unramified).  Then, by a result of Artin 
(\cite[Corollaire~10.4]{lmb}), the quotient sheaf $X/G$ in the
fppf topology is an algebraic space, and 
$G \times_{\spec k} X = X \times_{X/G} X$. 
There is a Zariski open dense subspace $V \subseteq X/G$ 
that is a scheme (\cite[Proposition~6.7]{knudson}); 
if $U$ is the inverse image of $V$ in $X$, then the restriction 
$U \to U/G = V$ is a $G$-torsor (i.e. a principal $G$-bundle 
in the fppf topology, cf{.} \cite{dg}). In the case where $X$ is
a $k$-variety, this is precisely the rational quotient map
we discussed in the introduction, i.e., $k(U/G) = k(X)^{G}$.
\end{remark}

\begin{prop} \label{no-name1}
Assume that $G$ is a group scheme of finite type over $k$, acting on a
quasi-separated $k$-scheme $X$, with a non-empty invariant open subscheme
on which the action is free.
Let $\cE$ be a $G$-equivariant locally 
free sheaf of rank $r$ on $X$. Then there exists a non-empty 
open $G$-invariant subscheme $U$ of $X$, such that 
the restriction $\cE\mid_{U}$ is isomorphic to the trivial 
$G$-equivariant sheaf $\mathcal{O}_{U}^{r}$.
\end{prop}

To see that Lemma~\ref{no-name} (over an arbitrary base field $k$)
follows from Proposition~\ref{no-name1}, 
recall the the well-known equivalence between the category 
of $G$-equivariant vector bundles on $X$ and the category 
of $G$-equivariant locally free sheaves on $X$. 
One passes from a $G$-bundle $V \to X$ to the $G$-equivariant 
locally free sheaf of sections of $V$; conversely,
to each $G$-equivariant locally free sheaf $\cE$ on $X$ one
associates the spectrum of the sheaf of symmetric algebras 
of the dual $\cE^{\vee}$ over $X$.

Note also that in the course of proving Lemma~\ref{no-name},
we may assume without loss of generality that $X$ is 
{\em primitive}, i.e., $G$ transitively permutes the irreducible 
components of $X$ (equivalently, $k(X)^G$ is a field). Indeed, an
arbitrary $G$-variety $X$ is easily seen to be birationally isomorphic 
to a disjoint union of primitive $G$-varieties $X_1, \dots, X_r$, and it 
suffices to prove Lemma~\ref{no-name} for each $X_i$. On the other hand,
if $X$ is primitive, then every non-empty $G$-invariant open subset
is dense in $X$. This shows that Lemma~\ref{no-name} follows 
from Proposition~\ref{no-name1}, as claimed. 

{\em Proof of Proposition~\ref{no-name1}.}
After replacing $X$ by a non-empty open subscheme
we may assume that the action of $G$ on $X$ is free. 
By passing to a dense invariant subscheme of $X$, we may assume that
$X/G$ is a scheme, and  $X \to X/G$ is a $G$-torsor; see 
Remark~\ref{restriction}.
By descent theory, the $G$-equivariant sheaf $\cE$ comes from 
a locally free sheaf $\cF$ on $X/G$;
see, for example,~\cite[Theorem~4.46]{descent}. 
By restricting to a non-empty open subscheme of $X/G$ once again,
we may assume that $\cF$ is isomorphic to $\mathcal{O}_{X/G}^{r}$. Then 
$\cE$ is $G$-equivariantly isomorphic to $\mathcal{O}_{X}^{r}$, as claimed.
\qed

\begin{remark} The same argument goes through if the base field $k$ 
(or, equivalently, the base scheme $\spec(k)$) is replaced by 
a algebraic space $B$, so that $X$ is defined over $B$, and the
group scheme $G$ is assumed to be flat and finitely presented 
over $B$.
\end{remark}

\section{Proof of Theorem~\ref{thm1}}

We will prove Theorem~\ref{thm1} in two steps: first in the case 
where $X = V$ is a generically free linear representation of $G$, 
then for arbitrary $X$.

\smallskip
{\bf Step 1:} Suppose $X = V$ is a generically free linear 
representation of $G$.

Let $m = \dim(V)$.  By the no-name lemma, 
there exist $G$-equivariant
birational isomorphisms $\alpha$ and $\beta$ such that the diagram
\[ 
\xymatrix{W \times \bbA^m \ar@{-->}[r]^{\alpha} \ar@{->}[d]^{pr_1} &
  W \times W^{\phi} \ar@{->}[d]^{pr_1} &
  W \times W^{\phi} \ar@{-->}[r]^{\beta} \ar@{->}[d]^{pr_2} &
  W^{\phi} \times \bbA^m \ar@{->}[d]^{pr_1}  \cr
  W \ar@{=}[r] \ar@{-->}[d]^{\pi_W} &
  W \ar@{-->}[d]^{\pi_W} &
  W^{\phi} \ar@{=}[r] \ar@{-->}[d]^{\pi_{W^{\phi}}} &
  W^{\phi} \ar@{-->}[d]^{\pi_{W^{\phi}}} \cr
   W/G \ar@{=}[r] &
   W/G \ar@{=}[r] &
   W^{\phi}/G \ar@{=}[r] &
   W^{\phi}/G}  \]
commutes. Now we can take 
$f = \beta \circ \alpha \colon 
V \times \bbA^m \dasharrow V^{\phi} \times \bbA^m$.

\smallskip
{\bf Step 2:} Suppose $X$ is an arbitrary generically free
$G$-variety. 

Let $V$ be a generically free linear representation of $G$ and
$p \colon X \times V \lra V$ be the projection onto the second 
factor. By the no-name lemma, $X \times V$ is
birationally isomorphic to $X \times \bbA^m$; this yields a dominant
rational map of $G$-varieties $X \times \bbA^m \dasharrow V$, which
we will continue to denote by $p$. After replacing $X$ by 
$X \times \bbA^m$, we may assume that there exists a dominant rational 
map $p \colon X \dasharrow V$. We now consider the commutative diagram
\[ 
\xymatrix{
  X \ar@{-->}[r]^{p} \ar@{-->}[d] & V \ar@{-->}[d]  \cr
  X/G \ar@{-->}[r]^{p/G}  & V/G,}
\]
where the vertical arrows are rational quotient maps. 
We claim that $X$ is birationally isomorphic to the fiber 
product $X/G \times_{V/G} V$,
where the $G$-action on this fiber product is induced from the $G$-action 
on $V$ (in other words, $G$ acts trivially on $X/G$ and on $V/G$). In the
case where $k$ is an algebraically closed field of characteristic zero,
this is proved in~\cite[Lemma 2.16]{reichstein}.
For general $k$, choose an open invariant subscheme $U$ of $V$ 
such that $G$ acts freely over $U$, the quotient $U/G$ is a scheme, 
and the projection $U \to U/G$ is a $G$-torsor; see Remark~\ref{restriction}. 
By restricting $X$, we may assume that the morphism $X \to X/G$ 
is also a $G$-torsor, and that $X$ maps into $U$. Then we get 
a commutative diagram
\[ 
\xymatrix{
  X \ar[r]^{p} \ar[d] & U \ar[d]  \cr
  X/G \ar[r]^{p/G}  & U/G}
\]
where the columns are $G$-torsors and the top row is $G$-equivariant. 
Any such diagram is well known to be cartesian; this proves our claim.

Similarly, $X^{\phi} \simeq X^{\phi}/G \times_{V^{\phi}/G} V^{\phi}$.
By Step 1 there is a $G$-equivariant birational isomorphism 
$f \colon V \times \bbA^m \dasharrow V^{\phi} \times \bbA^m$ which 
makes the diagram
\[ \xymatrix{
  & V \times \bbA^m \ar@{-->}[d]^{\pi_V} \ar@{-->}[dr]^f &    \cr
  X/G \ar@{-->}[r]^{p/G}  \ar@{-->}[dr]^{id} & V/G \ar@{-->}[dr]^{id}  & 
V^{\phi} \times \bbA^m \ar@{-->}[d]^{\pi_{V^{\phi}}}    \cr
  & X^{\phi}/G \ar@{-->}[r]^{p/G}  & V^{\phi}/G } 
\]
commute. Consequently, $f$ induces a $G$-equivariant birational isomorphism
between the fiber products $X/G  \times_{V/G} (V \times \bbA^m)$ and
$X^{\phi}/G  \times_{V/G} (V^{\phi} \times \bbA^m)$
i.e., between $X \times \bbA^m$ and $X^{\phi} \times \bbA^m$.
This completes the proof of Theorem~\ref{thm1}.

\begin{remark} \label{rem3.2} Our proof shows that the integer $m$ 
in the statement of Theorem~\ref{thm1} can be taken to 
be the minimal value of $2 \dim(W)$,
as $W$ ranges over the generically free linear representations of $G$. 
If $X$ is itself a generically free linear representation of $G$
then we can take $m$ to be the minimal value of $\dim(W)$, where
again $W$ ranges over the generically free linear 
representations of $G$; see the proof of Step 1.
\end{remark}

\section{Examples where $X$ and $X^{\phi}$ are birationally isomorphic}
\label{sect.so}

The question of whether or not
$X$ and $X^{\phi}$ are birationally isomorphic over $k$ is 
delicate in general. Birational isomorphism over $K = k(X)^G$
is more accessible because it can be restated in terms of 
Galois cohomology. In this section we will to show that 
in many cases $X$ and $X^{\phi}$ are, indeed, birationally 
isomorphic over $K$ (and thus over $k$).

Let $X$ be a primitive $G$-variety.  Recall that
$X$ is called primitive if $G$ transitively permutes the
irreducible components of $X$; see Section~\ref{sect.no-name}.
That is, $X$ is primitive
if $K = k(X)^G$ is a field or equivalently, if the rational quotient
variety $X/G$ is irreducible.
As we saw in Remark~\ref{restriction}, the rational quotient map 
$\pi \colon X \dasharrow X/G$ is a $G$-torsor over a non-empty
open subscheme of $X/G$; hence, the $G$-action on $X$ gives rise to
a Galois cohomology class in $H^1(K, G)$, which we shall denote by
$[X]$; cf., \cite{serre-gc}, \cite{popov}. 
 
An automorphism $\phi$ of $G$ induces an automorphism $\phi_*$ 
of the (pointed) cohomology set $H^1(K, G)$, where
$[X^{\phi}] = \phi_*([X])$. In particular, the $G$-varieties 
$X$ and $X^{\phi}$ are birationally isomorphic 
over $K$ if and only if $\phi_*([X]) = [X]$ in $H^1(K, G)$.

\begin{example} \label{ex3.14} If $[X] = 1$ then $\phi_*([X]) = 1$;
hence, $X$ and $X^{\phi}$ are birationally isomorphic. 
Explicitly,
in this case $X$ is birationally isomorphic to the ``split" $G$-variety
$Y \times G$, with $G$ acting on the second component by left translations,
and a birational isomorphism between $Y \times G$ and $(Y \times G)^{\phi}$
is given by $(y, g) \lra (y, \phi(g))$.

In particular, if $G$ is a special group, 
i.e., $H^1(K, G) = \{ 1 \}$ for every $K/k$, then
$X$ and $X^{\phi}$ are birationally isomorphic for every
generically free $G$-variety $X$. Examples of special groups 
are $G = \GL_n$, $\SL_n$, $\Sp_{2n}$; see~\cite[Chapter X]{serre-lf}.
\end{example}

The following lemma extends this simple argument a bit further.

\begin{lem} \label{lem3.12} Let $X$ be a primitive generically 
free $G$-variety and $K = k(X)^G$.
Suppose $[X]$ lies in the image of the natural map
$H^1(K, G_0) \lra H^1(K, G)$, where $G_0$ 
is a closed subgroup of $G$ such that 
$\phi_{| \, G_0} = \id \colon G_0 \lra G_0$.
Then $X$ and $X^{\phi}$ are birationally isomorphic as $K$-varieties. 
\end{lem}

\begin{proof} Let $i \colon G_0 \hookrightarrow G$ be 
the inclusion map.  The commutative diagram
\[ \xymatrix{
  G_0 \ar@{->}[r]^{i} \ar@{->}[d]^{id} & G \ar@{->}[d]^{\phi} \cr
  G_0 \ar@{->}[r]^{i} & G }
\]
of groups induces a commutative diagram of cohomology sets
\[ \xymatrix{
  H^1(K, G_0) \ar@{->}[r]^{i_*} \ar@{->}[d]^{id} & 
H^1(K, G) \ar@{->}[d]^{\phi_*} \cr
  H^1(K, G_0) \ar@{->}[r]^{i_*} & H^1(K, G). }
\]
Since $[X]$ is in the image of $i_*$, this diagram shows that
$\phi_*([X]) = [X]$.
\end{proof}

Note that if $[X] = 1$ then $[X]$ is the image 
of the trivial element of $H^1(K, G_0)$,
where $G_0 = \{ 1 \}$.  Example~\ref{ex3.14} is thus
a special case of Lemma~\ref{lem3.12}. We now turn to a more 
sophisticated application of Lemma~\ref{lem3.12} (with 
non-trivial $G_0$). 

\begin{prop} \label{prop.so}
Let $G$ be the special orthogonal group $\SO(q)$, where $q$
is a non-degenerate isotropic $n$-dimensional quadratic form 
defined over $k$, and let $X$ be an irreducible generically 
free $G$-variety.  Assume $n \ne 4$ and $\Char(k) \ne 2$.
Then $X$ and $X^{\phi}$ are birationally isomorphic 
over $K = k(X)^{\SO(q)}$ (and hence, over $k$)
for any automorphism $\phi$ of $\SO(q)$.
\end{prop}

\begin{proof} As we remarked in the introduction, we are allowed
to replace $\phi$ by $\phi' = \phi \circ \inn_h$, where $\inn_h$ denotes
conjugation by $h \in \SO(q)$; indeed, $X^{\phi}$ and $X^{\phi'}$
are isomorphic (over $K$) via $x \mapsto h \cdot x$. In particular, 
the proposition is true if $\phi$ is an inner automorphism of $\SO(q)$.

It is well known that every automorphism $\phi \colon \SO(q) \lra \SO(q)$ 
has the form $g \lra h g h^{-1}$, where $h \in \Orth(q)$; cf., e.g.,
\cite[Section XI]{dieudonne}.  If $n$ is odd then, after replacing $h$ 
by $\det(h) h$, we may assume that $h \in \SO(q)$,
i.e., $\phi$ is inner. This proves the proposition for odd $n$.

 From now on we will assume that $n \ge 6$ is even and $\det(h) = -1$.
We may also assume that
$q(x_1, \dots, x_n) = a_1 x_{1}^{2} + \dots + a_n x_{n}^{2}$ for
some $a_1, \dots, a_n \in k^*$ and, after 
composing $\phi$ with an inner automorphism of $\SO(q)$,
that $h = \diag(-1, 1, \dots, q)$.
Let $D_0 \simeq (\bbZ/2 \bbZ)^{n-1}$, $D \simeq (\bbZ/2 \bbZ)^n$
be the subgroups of diagonal 
matrices in $\SO(q)$, $\Orth(q)$ respectively, and
$i \colon D_0 \hookrightarrow \SO(q)$,
$j \colon D \hookrightarrow \SO(q)$
be the natural inclusion maps. Note that $\phi$ restricts to 
a trivial automorphism of $D_0$. Thus, in view of Lemma~\ref{lem3.12}, 
it suffices to show that $i_* \colon H^1(K, D_0) \lra H^1(K, \SO(q))$ 
is surjective.

Consider the commutative diagram
\[ \xymatrix{ 1 \ar@{->}[r] & D_0 \ar@{->}[r] \ar@{->}[d]^i & 
D \ar@{->}[d]^j \ar@{->}[r]^{\det} & \bbZ/2 \bbZ \ar@{->}[r] & 1 \cr
  & \SO(q) \ar@{->}[r] & 
\Orth(q) &  & } 
\] 
of algebraic groups and the induced commutative diagram 
\[
\xymatrix{ 1 \ar@{->}[r] & H^1(K, D_0) \ar@{->}[r] \ar@{->}[d]^{i_*} & 
(K^*/(K^*)^2)^n \ar@{->}[d]^{j_*} \ar@{->}[r]^{p} & K^*/(K^*)^2
\ar@{->}[r] & 1 \cr
  & H^1(K, \SO(q)) \ar@{->}[r] & 
H^1(K, \Orth(q)) &  & } 
\]
in cohomology. The top row is an exact sequence of abelian groups;
$H^1(K, D_0)$ can thus be identified with the kernel of the product map
$p \colon (K^*/(K^*)^2)^n \lra K^*/(K^*)^2$, where $p(b_1, \dots, b_n)
= b_1 \dots b_n \pmod{(K^*)^2}$.

Now recall that
$H^1(K, \Orth(q))$ is in a natural 1-1 correspondence with
isometry classes of $n$-dimensional quadratic forms $q'$ and that
$j_*$ takes $(b_{1}, \dots, b_n) \in (K^{*}/{K^{*}}^{2})^n$ to
the quadratic form $q' = a_1 b_1 x_1^2 + \dots + a_n b_n x_n^2$.
Similarly, $H^1(K, \SO(q))$ is in a natural 1-1 correspondence with
isometry classes of $n$-dimensional quadratic forms $q'$ such
that $q'$ has the same discriminant as $q$, and $i_*$ takes
$(b_{1}, \dots, b_n) \in (K^{*}/{K^{*}}^{2})^n$,  
with $b_1 \dots b_n = 1$ in $K^{*}/{K^{*}}^{2}$,
to $q' = a_1 b_1 x_1^2 + \dots + a_n b_n x_n^2$;
cf., e.g., \cite[III, Appendix 2, \S 2]{serre-gc}. It is 
clear from this description that both $i_*$ and $j_*$ are
surjective.
\end{proof}

\section{Proof of Theorem~\ref{thm2}}

Recall that elements of $H^1(K, \PGLn)$ are in a natural 1-1 correspondence 
with

\smallskip
(i) generically free $\PGLn$-varieties $X$, with $k(X)^{\PGLn} = K$,
up to birational isomorphism over $K$, or alternatively, with

\smallskip
(ii) central simple algebras $A/K$ of degree $n$, up to $K$-isomorphism;

\smallskip
\noindent
see~\cite{serre-gc}, \cite{popov}, \cite{rv}.
We will denote the central simple 
algebra corresponding
to an irreducible generically free $\PGLn$-variety $X$ (respectively, to an
element $\alpha \in H^1(K, \PGLn)$ by $A_X$ (respectively, by $A_{\alpha}$). 
If $\phi \colon \PGLn \lra \PGLn$ is the automorphism
given by $g \lra (g^{-1})^{\rm transpose}$ then 
$A_{\phi_*(\alpha)}$ is the opposite algebra $A_{\alpha}^{op}$;
cf. e.g., \cite[pp. 152-153]{serre-lf}.  In other words, 
$A_{X^{\phi}} = A_X^{op}$.  The following lemma gives a necessary 
and sufficient conditions
for $X$ and $X^{\phi}$ to be birationally isomorphic over $K$.

\begin{lem} \label{lem.exp2} 
Let $X$ be an irreducible generically free $\PGLn$-variety.
Then the following conditions are equivalent.

\smallskip
(a) $X$ and $X^{\phi}$ are birationally isomorphic over $K = k(X)^{\PGLn}$,

\smallskip
(b) $A_X$ is $K$-isomorphic to $A_X^{op}$,

\smallskip
(c) $A_X$ has exponent $1$ or $2$ in the Brauer group $\Br(K)$.
\end{lem}

\begin{proof} (a) $\Longleftrightarrow$ (b):
Let $\alpha = [X] \in H^1(K, \PGLn)$. Then
as we observed in Section~\ref{sect.so}, 
$[X^{\phi}] = \phi_{\ast}(\alpha)$. Thus 
$A_{X^{\phi}} = A_{\phi_*(\alpha)} = A_{\alpha}^{op} = A_X^{op}$,
so that $X$ and $X^{\phi}$ are
birationally isomorphic over $K$ if and only if $A_X$ is $K$-isomorphic
to $A_X^{op}$.

The equivalence of (b) and (c) is obvious, since $A^{op}$ is the inverse
of $A$ in $\Br(K)$.
\end{proof}

Lemma~\ref{lem.exp2} does not directly address the question we
are interested in, namely, the question of whether $X$ 
and $X^{\phi}$ are birationally isomorphic over the base field $k$.
Note however, that a birational isomorphism $\alpha \colon X \lra X^{\phi}$ 
defined over $k$, restricts to a $k$-automorphism of the field of
invariants $K = k(X)^{\PGLn} = k(X^{\phi})^{\PGLn}$. 
Our proof of Theorem~\ref{thm2} is based on
the observation that if $\Aut_k(K) = \{ 1 \}$, then

$X$ and $X^{\phi}$ are isomorphic over $k$ $\Longleftrightarrow$

$X$ and $X^{\phi}$ are isomorphic over $K$ $\Longleftrightarrow$

$A_X$ has exponent $1$ or $2$ in the Brauer group of $K$.

\smallskip
Thus in order to prove Theorem~\ref{thm2} it suffices to
construct (i) a finitely generated field extension $K/k$
such that $\Aut_k(K) = \{ 1 \}$ and (ii) a central simple
algebra $A/K$ of degree $n$ and exponent $n$. These
constructions are carried out in Lemmas~\ref{lem3.3}
and~\ref{lem3.5} below. To simplify the exposition, we will
state Lemmas~\ref{lem3.3} and~\ref{lem3.5} for an
algebraically closed ground field $k$.  
We will then explain how to modify our construction to
make it work over any infinite field $k$ containing 
a primitive $n$th root of unity.

\begin{lem} \label{lem3.3} 
Suppose $k$ is an algebraically closed field. Then
there exists an algebraic surface $S/k$ which admits 
no non-trivial birational automorphisms. In other words,
$\Aut_k \, k(S) = \{ 1 \}$.
\end{lem}

\begin{proof} Choose two smooth non-isomorphic curves $C_1$ 
and $C_2$, of genus $g_1$ and $g_2$ respectively (say, $g_1 \ge g_2$), 
with no non-trivial automorphisms, and set $S = C_1 \times C_2$.

We claim that every birational automorphism $f \colon S \dasharrow S$ 
is trivial. To prove this, note that $f$ restricts to 
a regular map $C \lra S$ for every smooth 
curve $C$ in $S$. Taking $C = C_1 \times \{ y \}$ for
some $y \in C_2$, we see that $pr_2 \circ f$ is a regular map $C \lra C_2$.
(Here $pr_2 \colon S \lra C_2$ is the projection to the second factor.)
Since $g_1 \ge g_2$, the Hurwitz formula tells us that this map cannot 
be dominant, i.e., it sends $C$ to a single point. In other words,
$f (C_1 \times \{ y \}) \subset C_1 \times \{ y' \}$ for 
some $y' \in C_2$. Since $f$ is a birational automorphism,
$f (C_1 \times \{ y \})$ can be a single point for only finitely many 
$y \in C_2$. For every other $y \in C_2$ there exists a $y' \in C_2$
such that
\[ f (C_1 \times \{ y \}) = C_1 \times \{ y' \} \, . \]
Applying the Hurwitz formula once again, we see that
$f$ induces an isomorphism between
$C_1 \times \{ y \}$ and $C_1 \times \{ y' \}$. Since $C_1$
has no non-trivial automorphisms, this isomorphism 
is given by $f(x, y) = (x, y')$.
Equivalently, for $x \in C_1$, $f$ restricts to a morphism 
$\{ x \} \times C_2 \lra \{ x \} \times C_2$. The Hurwitz formula 
now tells us that this map is an automorphism. Since $C_2$ has no
non-trivial automorphisms, we conclude that $f(x, y) = (x, y)$, 
for every $x, y \in S$.
\end{proof}

\begin{lem} \label{lem3.5}
Assume $k$ is an algebraically closed field, $\Char(k)$ 
does not divide $n$ and $K/k$ be a finitely generated field extension
of transcendence degree $\ge 2$. Then for every $n \ge 3$
there exists a division algebra $D/K$ of degree $n$ and exponent $n$.
\end{lem}

\begin{proof}
Consider a model $X$ for $K$, i.e., an algebraic variety $X$ with function
field $k(X) = K$. Choose a smooth point $x \in X$ and
a system of local parameters $t_1, \dots, t_d$ in the local ring $\cO$;
here $d = \dim(X) = \trdeg_k(K) \ge 2$.
We claim that the symbol algebra
$D = (t_1, t_2)_n$, i.e., the $K$-algebra given by
generators $x_1, x_2$ and relations $x_1^n =t_1$, $x_2^n = t_2$, 
$x_1 x_2 = \zeta x_2 x_1$ has exponent $n$ in $\Br(K)$.
(Here $\zeta_n$ is a primitive $n$th root of unity in $k$.)

To prove this, consider the completion $\ccO = k[[t_1, \dots, t_d]]$
of the local ring $\cO$, where $k[[t_1, \dots, t_d]]$ denotes
the ring of formal power series in the variables $t_1, \dots, t_d$.
Note that $\cO \subset \ccO$ and thus, after passing to
the fields of fractions, $K \subset k((t_1, \dots, t_d))$.
The image of $D$ under the restriction
map $\Br(K) \lra \Br \big( k((t_1, \dots, t_d)) \big)$ is the symbol
algebra $D' = (t_1, t_2)_n$ over $k((t_1, \dots, t_d))$. A simple
valuation-theoretic argument shows that $D'$ has exponent $n$;
cf.~\cite[Proposition 3.3.26]{rowen}.
Hence, so does $D$.
\end{proof}

Lemmas~\ref{lem3.3} and~\ref{lem3.5} complete the proof 
of Theorem~\ref{thm2} in the case where the base 
field $k$ is algebraically closed and $\Char(k)$ does not divide $n$. 
Then the same argument will work over a non-closed field 
$k$ with a primitive $n$th root of unity, if
we can choose the curves $C_i$ ($i = 1$, $2$) 
in Lemma~\ref{lem3.3} so that $C_i$ each is defined over $k$ 
and has a $k$-point $p_i$. (The primitive $n$th root of unity is 
needed to define the symbol algebra $D$ in the proof of Lemma~\ref{lem3.5}.)
Taking the $k$-variety $X$ in Lemma~\ref{lem3.5} to be $C_1 \times C_2$,
we see that $x = (p_1, p_2)$ is a smooth $k$-point of $X$, and
the rest of the proof of Lemma~\ref{lem3.5} goes through unchanged.

To construct the curve $C_i$ as above (for $i = 1, 2$),
fix a $k$-point $p$ in $\bbP^2$ and $C_i$ to be the general element of the
$k$-linear system of degree $d_i$ curves passing through $p_i$. If $k$
is an infinite field then $C_i$ is a smooth curve of 
genus $G_i = \frac{1}{2} (d_i -1) (d_i - 2)$. By construction,
$C$ has a $k$-point $p_i$. Moreover, if $d_i \ge 4$, then $C_i$ has 
no non-trivial automorphisms.

This completes the proof of Theorem~\ref{thm2}.
\qed

\begin{remark} \label{rem.m}
Theorem~\ref{thm2} remains valid if the condition that 
$k$ contains a primitive $n$th root of unity is replaced by
the (weaker) condition that $k$ contains a primitive $m$th root of unity 
for some divisor $m$ of $n$ such that $m \ge 3$. The proof is the same, 
except that instead of the symbol algebra $D = (t_1, t_2)_n$ of degree $n$
and exponent $n$, we use the algebra $\Mat_{n/m}(E)$ of degree $n$
and exponent $m$, where $E = (t_1, t_2)_m$.
\end{remark}
 
\begin{remark} \label{n=2} 
Theorem~\ref{thm2} fails for $n = 2$; indeed, $A$ and $A^{op}$ 
are isomorphic over $K$ for any central simple 
algebra $A/K$ of degree $2$. Alternatively, 
$\PGL_{2} \simeq \SO_{3}$, so if Theorem~\ref{thm2} were 
true for $n = 2$, it would contradict Proposition~\ref{prop.so}.
\end{remark}

\section{Further examples}

In this section we will assume that $G$ is a finite group and
$k$ is an algebraically closed field of characteristic zero.

\begin{prop} \label{prop1} 
(a) For every finitely generated field extension $K/k$ and every
finite group $G$ there exists a $G$-Galois extension $L/K$.
Equivalently, there exists an irreducible
$G$-variety $X$ such that $k(X)^G = K$.

(b) Suppose $\Aut_k(K) = \{ 1 \}$ and $\phi \colon G \lra G$ is an outer
automorphism of a finite group $G$. Then for every $X$, as in (a),
the $G$-varieties $X$ and $X^{\phi}$ are not birationally isomorphic. 
\end{prop}

\begin{proof} (a) By the Riemann existence theorem there exists
a $G$-Galois extension $L_0/k(t)$, where $t$ is an independent variable.
Hence, there exists a $G$-Galois extension $L_1/K(t)$, 
where $L_1 = L \otimes_{k(t)} K(t)$.  The Hilbert irreducibility 
theorem now allows to construct a $G$-Galois 
extension $L/K$ by suitably specializing $t$ in $K$.

(b) 
Irreducible $G$-varieties $X$ (up to birational isomorphism)
such that $k(X)^G = K$, are in 1-1 correspondence with 
$G$-Galois field extensions $L/K$.
A birational isomorphism $\alpha \colon X \dasharrow X^{\phi}$
of $G$-varieties induces an isomorphism 
\[ \xymatrix{
  L = k(X) \ar@{->}[r]^{\alpha} \ar@{-}[d] & k(X^{\phi}) = L \ar@{-}[d] \cr
  K \ar@{=}[r] & K.} 
\]
Then $\alpha \in \Gal(L/K) = G$, and since the above diagram commutes,
we have $\alpha g(l) = \phi(g) \alpha(l)$ for every $g \in G$ and $l \in L$.
In other words, $\phi(g) = \alpha g \alpha^{-1}$, contradicting our
assumption that $\phi$ is an outer automorphism.
\end{proof}

\smallskip
\noindent{\bf Acknowledgement.}
The first-named author would like to thank V. L. Popov for
helpful discussions.

\end{document}